\documentclass{amsart}
\title[Combinatorial Identities Via Phi Functions...]{Combinatorial Identities Via Phi Functions and Relatively Prime Subsets}
\usepackage{amssymb,amsmath,amsthm,epsfig,graphics,latexsym}
 \theoremstyle{definition}
 \newtheorem{definition}{Definition}
  \theoremstyle{plain}

  \newtheorem{theorem}    {Theorem}
  \newtheorem{corollary}  {Corollary}

  \theoremstyle{remark}

\begin{document}
  \author{Mohamed El Bachraoui}
  \address{Dept. Math. Sci,
 United Arab Emirates University, PO Box 17551, Al-Ain, UAE}
 \email{melbachraoui@uaeu.ac.ae}
 \keywords{Mertens function, M\"obius function, Combinatorial identities, Phi functions, Relatively prime sets}
 \subjclass{11A25, 11B05, 11B75}
  %
  \begin{abstract}
  Let $n$ be a positive integer and let $A$ be nonempty finite set of positive integers. We say that $A$ is relatively prime if
  $\gcd(A) =1$ and that $A$ is relatively prime to $n$ if $\gcd(A,n)=1$. In this work we count 
  the number of nonempty subsets of $A$ which are relatively prime and the number of nonempty subsets of
  $A$ which are relatively prime to $n$. Related formulas are also obtained for the number of such subsets having
  some fixed cardinality. This extends previous work for the cases where $A$ is an interval or a set
  in arithmetic progression. Applications include:

  \noindent
  a) An exact formula is obtained for the number of elements of $A$ which are co-prime to $n$; note that this number is $\phi(n)$
  if $A=[1,n]$.

  \noindent
  b) Algebraic characterizations are found for a nonempty finite set of positive integers to have elements which are all
   pairwise co-prime and consequently a formula is given for the number of nonempty subsets of $A$ whose elements are pairwise co-prime.

  \noindent
  c) We provide combinatorial formulas involving Mertens function.
  \end{abstract}
  \date{\textit{\today}}
  \maketitle
 \section{Introduction}
 Throughout let $n$ and $\alpha$ be positive integers and let $A$ be a nonempty finite set of positive integers.
 Let $\mu$ be the M\"obius function and let
 $\lfloor x \rfloor$ be the floor of $x$. If $a$ and $b$ are positive integers such that $a\leq b$, then we let $[a,b]=\{a,a+1,\ldots,b\}$.
 Further we need the following set theoretical notation. Let $\# X = |X|$ denote the cardinality of a set $X$, let $\mathcal{P}(X)$ denote the power set of $X$, let $\mathcal{P}^{\ast}(X)$ denote
 the set of nonempty subsets of $X$, and let
 $\mathcal{P}_{\alpha}(X)$ denote the set of subsets of $X$ whose cardinality is $\alpha$.
 Multiples of integers and their cardinality are crucial in this paper.
 \begin{definition}
 For any positive integer $d$ let $V(A,d)$ be the set of multiples of $d$ in $A$ and let
 $v(A,d)=|V(A,d)|$.
 \end{definition}
 \begin{theorem} \label{multiples}
 We have
 \[ v(A,d) = \sum_{a\in A} (\lfloor a/d \rfloor - \lfloor (a-1)/d \rfloor ). \]
 \end{theorem}
 \begin{proof}
 The result follows since
 \[
 \lfloor a/d \rfloor - \lfloor (a-1)/d \rfloor =
  \begin{cases}
  1,\ \text{if $a \in V(A,d)$} \\
  0,\ \text{if $a \not\in V(A,d)$.}
  \end{cases}
  \]
 \end{proof}
  We now state some of the results which we shall prove in this work. \\
  \noindent
 1) While the number of elements in the set $[1,n]$ which are co-prime to $n$ is $\phi(n)$, to the author's
  knowledge no such formula exists for the number of elements in an arbitrary nonempty finite set $A$ of positive integers.
  In this paper we will show in Corollary \ref{coprime-n} that such a number is
 \[
  \sum_{d|n}\mu(d) v(A,d).
 \]
 \noindent
 2) Cameron and Erd\H{o}s in \cite{Cameron-Erdos} considered for any positive real number $x$ the sets of positive integers
 whose elements are $\leq x$ and pairwise co-prime and considered the related sets of positive integers with elements $\leq x$ and which are
 free of co-prime pairs.
 The authors gave asymptotic lower and upper bounds for the numbers of these sets. Calkin and Granville
 in \cite{Calkin-Granville} gave asymptotic formulas for such numbers and improved the result of Cameron and Erd\H{o}s. In this paper
 we shall prove in Corollary \ref{cor:pairwise-coprime} that the number of subsets of $A$ whose elements are pairwise co-prime is:
 \[ \# \left\{B\subseteq A:\ B\not=\emptyset\ \text{and\ }
   2|B| -1 = \left( 1 + 4 \sum_{d=1}^{\sup B} \mu(d) v(B,d) (v(B,d)-1) \right)^{1/2} \right\}.
 \]
 \noindent
 3) From Corollary \ref{subset-relprime}, for any nonempty subset
 $B$ of $A$ satisfying $\alpha \leq |B|$ the identity
 \[ \binom{|A|}{\alpha}=\sum_{d=1}^{\sup A}\mu(d)\binom{v(A,d)}{\alpha}, \]
 yields
 \[ \binom{|B|}{\alpha}=\sum_{d=1}^{\sup B}\mu(d)\binom{v(B,d)}{\alpha}. \]
 As a consequence, for any nonempty subset $B$ of $A$, if
 \[ 1 + 4 \sum_{d=1}^{\sup A} \mu(d) v(A,d)\left( v(A,d)-1 \right) \]
 is a square, then so is
 \[ 1 + 4 \sum_{d=1}^{\sup B} \mu(d) v(B,d)\left( v(B,d)-1 \right). \]
 \noindent
 4) In Theorem \ref{Mertens-m,n} we will show that for any positive integers $1<m\leq n$
 \[  \sum_{d=1}^n \mu(d)
  (2^{\lfloor\frac{n}{d} \rfloor - \lfloor \frac{n-1}{d} \rfloor + \lfloor\frac{m}{d} \rfloor - \lfloor \frac{m-1}{d} \rfloor} )=
  \begin{cases}
 M(n)\ \text{if $(m,n) >1$,} \\
 1 + M(n)\ \text{if $(m,n)=1$,}
 \end{cases}
 \]
 where $M$ is Mertens function given by
 $M(n) =\sum_{d=1}^{n} \mu(d)$.

 \noindent
 5) A direct consequence of Theorem \ref{Mertens-for-composite} is the following combinatorial identity.
 For simplicity of notation if $b$ is a positive integer, then
 \[ b A = \{b a:\ a \in A \}. \]
 Let $c$ be a composite positive integer, say $c=ab$ with $a >1$ and $b >1$. Then
 for any finite sets $A=\{a_1,a_2,\ldots,a_k \}$ and $B=\{b_1,b_2,\ldots,b_l \}$ of positive integers with
 $\sup A = a$ and $\sup B = b$ we have
 \[
 M(c) = \sum_{d=1}^c \mu(d) 2^{v(b A,d)} =
  \sum_{d=1}^c \mu(d) 2^{v(\{a B,d)}.
  \]
 \noindent
  Our main tools are relatively prime sets and phi functions for sets of positive integers.
  The set $A$ is called \emph{relatively prime} if $\gcd(A) = 1$ and it is called
 \emph{relatively prime to $n$} if $\gcd(A\cup \{n\}) = \gcd(A,n) = 1$.
 From now on we assume that $\alpha \leq |A|$.
 \begin{definition}
  Let
  \[
  \begin{split}
  f(A) &= \# \{X\subseteq A:\ X\not=
  \emptyset\ \text{and\ } \gcd(X) = 1 \}, \\
  f_{\alpha} (A) &= \# \{X\subseteq A:\ \# X=
  \alpha \ \text{and\ } \gcd(X) = 1 \}, \\
  \Phi(A,n) &= \# \{X\subseteq A:\ X\not=
  \emptyset\ \text{and\ } \gcd(X,n) = 1 \}, \\
  \Phi_{\alpha} (A,n) &= \# \{X\subseteq A:\ \# X=
  \alpha \ \text{and\ } \gcd(X,n) = 1 \}.
  \end{split}
  \]
 \end{definition}
 Nathanson in \cite{Nathanson} introduced among others the functions $f(n)$, $f_{\alpha}(n)$, $\Phi(n)$, and $\Phi_{\alpha}(n)$
 (in our terminology $f([1,n])$, $f_{\alpha}([1,n])$, $\Phi([1,n],n)$, and $\Phi_{\alpha}([1,n],n)$ respectively) and found exact formulas
 along with asymptotic estimates for each of these functions.
 Formulas for these functions are found in \cite{ElBachraoui1, Nathanson-Orosz} for $A= [m,n]$
  and in \cite{ElBachraoui2} for $A=[1,m]$.
 Ayad and Kihel in \cite{Ayad-Kihel1, Ayad-Kihel2} considered extensions to sets in arithmetic progression and obtained
 identities for these functions for $A=[l,m])$ as consequences.
 Recently in \cite{ElBachraoui3, ElBachraoui-Salim} these functions have been studied for the union of two intervals,
 the special case of $A=[l,m]$ has been obtained, and
 various combinatorial identities involving M\"{o}bius and Mertens functions have been found as applications.
 For the purpose of this work we give these functions for $A = [l,m]$.
 \begin{theorem} \label{relprime-interval} 
 We have
 \[
 \begin{split}
 (a)\quad f([l,m]) &= \sum_{d=1}^{m}(2^{\lfloor\frac{m}{d}\rfloor - \lfloor\frac{l-1}{d}\rfloor}-1),\\
 (b)\quad f_{\alpha}([l,m]) &= \sum_{d=1}^{m}\binom{\lfloor\frac{m}{d}\rfloor - \lfloor\frac{l-1}{d}\rfloor}{\alpha}, \\
 (c)\quad \Phi([l,m],n) &= \sum_{d|n} 2^{\lfloor\frac{m}{d}\rfloor - \lfloor\frac{l-1}{d}\rfloor}, \\
 (d)\quad \phi_{\alpha}([l,m],n) &= \sum_{d|n} \binom{\lfloor\frac{m}{d}\rfloor - \lfloor\frac{l-1}{d}\rfloor}{\alpha}.
 \end{split}
 \]
 \end{theorem}
 An analysis of the functions $f$, $f_{\alpha}$, $\Phi$, and $\Phi_{\alpha}$ obtained for different
 cases of the set $A$
 lead us to more general formulas for any nonempty finite set of positive integers. See Sections
 Section \ref{sec:function-phi} and \ref{sec:function-f} below.
  \section{Phi functions for integer sets} \label{sec:function-phi}
 \begin{theorem} \label{main1}
 We have
 \[ (a)\quad
  \Phi(A,n)= \sum_{d|n}\mu(d) 2^{v(A,d)} = \sum_{d|n} \mu(d) |\mathcal{P}\left( V(A,d) \right)|,
 \]
 \[ (b)\quad
  \Phi_{\alpha}(A,n)= \sum_{d|n}\mu(d) \binom{v(A,d)}{\alpha} = \sum_{d|n}\mu(d)| \mathcal{P}_{\alpha} \left( V(A,d) \right)|.
 \]
 \end{theorem}
 \begin{proof}
  The second identities in (a) and (b) are trivial. As to the first identities we use
  induction on $|A|$. If $A=\{a \}=[a,a]$, then by Theorem \ref{relprime-interval} (c, d)
  \[ \Phi(A,n) =\sum_{d|n} \mu(d) 2^{\lfloor \frac{a}{d} \rfloor - \lfloor \frac{a-1}{d} \rfloor}\
  \text{and\ }
  \Phi_{\alpha}(A,n) = \sum_{d|n} \mu(d) \binom{\lfloor \frac{a}{d} \rfloor - \lfloor \frac{a-1}{d} \rfloor}{\alpha}. \]
  Assume that $A=\{a_1,a_2,\ldots,a_k \}$ and that the identities hold for $\{a_2,\ldots,a_k\}$. Then as to (a) we have
  \[
  \begin{split}
  \Phi(\{a_1,\ldots,a_k\},n)
  &=
  \Phi(\{a_2,\ldots,a_k\},n) + \Phi(\{a_2,\ldots,a_k\},\gcd(a_1,n) ) \\
  &=
  \sum_{d|n}\mu(d) 2^{\sum_{i=2}^k (\lfloor \frac{a_i}{d}\rfloor - \lfloor \frac{a_i -1}{d}\rfloor)}
   +
   \sum_{d|(a_1,n)}\mu(d) 2^{\sum_{i=2}^k (\lfloor \frac{a_i}{d}\rfloor - \lfloor \frac{a_i -1}{d}\rfloor)} \\
  &=
  2 \sum_{d|(a_1,n)}\mu(d) 2^{\sum_{i=2}^k (\lfloor \frac{a_i}{d}\rfloor - \lfloor \frac{a_i -1}{d}\rfloor)}
   +
   \sum_{\substack{d|n\\ d\nmid a_1}}\mu(d) 2^{\sum_{i=2}^k (\lfloor \frac{a_i}{d}\rfloor - \lfloor \frac{a_i -1}{d}\rfloor)} \\
  &=
  \sum_{d|(a_1,n)}\mu(d) 2^{\lfloor \frac{a_1}{d} \rfloor - \lfloor \frac{a_1-1}{d} \rfloor}
  2^{\sum_{i=2}^k (\lfloor \frac{a_i}{d}\rfloor - \lfloor \frac{a_i -1}{d}\rfloor)}
   +
   \sum_{\substack{d|n\\ d\nmid a_1}}\mu(d) 2^{\sum_{i=1}^k (\lfloor \frac{a_i}{d}\rfloor - \lfloor \frac{a_i -1}{d}\rfloor)} \\
  &=
  \sum_{d|(a_1,n)}\mu(d) 2^{\sum_{i=1}^k (\lfloor \frac{a_i}{d}\rfloor - \lfloor \frac{a_i -1}{d}\rfloor)}
   +
   \sum_{\substack{d|n\\ d\nmid a_1}}\mu(d) 2^{\sum_{i=1}^k (\lfloor \frac{a_i}{d}\rfloor - \lfloor \frac{a_i -1}{d}\rfloor)} \\
  &=
  \sum_{d|n}\mu(d) 2^{\sum_{i=1}^k
   (\lfloor \frac{a_i}{d}\rfloor - \lfloor \frac{a_i -1}{d}\rfloor)}.
   \end{split}
   \]
   As to (b) we have
   \[
   \begin{split}
   \Phi_{\alpha}(\{a_1,\ldots,a_k\},n)
  &=
  \Phi_{\alpha}(\{a_2,\ldots,a_k\},n) + \Phi_{\alpha -1}(\{a_2,\ldots,a_k\},\gcd(a_1,n) ) \\
  &=
  \sum_{d|n}\mu(d) \binom{\sum_{i=2}^k (\lfloor \frac{a_i}{d}\rfloor - \lfloor \frac{a_i -1}{d}\rfloor)}{\alpha}
   +
  \sum_{d|(a_1,n)}\mu(d) \binom{\sum_{i=2}^k (\lfloor \frac{a_i}{d}\rfloor - \lfloor \frac{a_i -1}{d}\rfloor)}{\alpha -1} \\
  &=
  \sum_{d|(a_1,n)}\mu(d)\left( \binom{\sum_{i=2}^k (\lfloor \frac{a_i}{d}\rfloor - \lfloor \frac{a_i -1}{d}\rfloor)}{\alpha}
  + \binom{\sum_{i=2}^k (\lfloor \frac{a_i}{d}\rfloor - \lfloor \frac{a_i -1}{d}\rfloor)}{\alpha -1} \right) + \\
  & \quad\quad
   \sum_{\substack{d|n\\ d\nmid a_1}}\mu(d) \binom{\sum_{i=2}^k (\lfloor \frac{a_i}{d}\rfloor - \lfloor \frac{a_i -1}{d}\rfloor)}{\alpha} \\
  &=
  \sum_{d|(a_1,n)}\mu(d) \binom{1+ \sum_{i=2}^k (\lfloor \frac{a_i}{d}\rfloor - \lfloor \frac{a_i -1}{d}\rfloor)}{\alpha}
  +
  \sum_{\substack{d|n\\ d\nmid a_1}}\mu(d) \binom{\sum_{i=2}^k (\lfloor \frac{a_i}{d}\rfloor - \lfloor \frac{a_i -1}{d}\rfloor)}{\alpha} \\
  &=
  \sum_{d|(a_1,n)}\mu(d) \binom{\sum_{i=1}^k (\lfloor \frac{a_i}{d}\rfloor - \lfloor \frac{a_i -1}{d}\rfloor)}{\alpha}
  +
  \sum_{\substack{d|n\\ d\nmid a_1}}\mu(d) \binom{\sum_{i=1}^k (\lfloor \frac{a_i}{d}\rfloor - \lfloor \frac{a_i -1}{d}\rfloor)}{\alpha} \\
  &=
  \sum_{d|n}\mu(d) \binom{\sum_{i=1}^k (\lfloor \frac{a_i}{d}\rfloor - \lfloor \frac{a_i -1}{d}\rfloor)}{\alpha}.
  \end{split}
  \]
   This completes the proof.
   \end{proof}
 \begin{corollary} \label{coprime-n}
  The number of positive integers in the set $A$
  which are co-prime to $n$ is
 \[
  \sum_{d|n}\mu(d) v(A,d).
 \]
 \end{corollary}
 \begin{proof}
 Apply Theorem \ref{main1}(b) with $\alpha =1$.
 \end{proof}
 \begin{corollary} \label{parity}
 If $n \in A$, then $\Phi(A,n) \equiv 0 \bmod 2$.
 \end{corollary}
 \begin{proof}
 If $n\in A$, then evidently $v(A,d) >0$ for all divisor $d$ of $n$ and thus the required congruence follows
 by Theorem \ref{main1}(a).
 \end{proof}
 \section{Relatively prime subsets of integer sets} \label{sec:function-f}
 \begin{theorem} \label{main2}
  We have
 \[ (a)\quad
  f(A)= \sum_{d=1}^{\sup A}\mu(d) \left(2^{v(A,d)} -1\right) = \sum_{d=1}^{\sup A}\mu(d) |\mathcal{P}^{\ast}\left( V(A,d)\right)|,
 \]
 \[ (b)\quad
  f_{\alpha}(A)= \sum_{d=1}^{\sup A}\mu(d) \binom{v(A,d)}{\alpha} =\sum_{d=1}^{\sup A}\mu(d) | \mathcal{P}_{\alpha}\left( V(A,d)\right) | .
 \]
 \end{theorem}
 \begin{proof}
 The second equalities in (a) and (b) are evident. As to the first identities, we use induction on $|A|$.
 If $A=\{a \}=[a,a]$, then by Theorem \ref{relprime-interval} (a, b)
  \[ f(A)= \sum_{d=1}^a \mu(d) \left(2^{\lfloor \frac{a}{d} \rfloor - \lfloor \frac{a-1}{d} \rfloor}-1\right)\ \text{and\ }
  f_{\alpha}(A) = \sum_{d=1}^a \mu(d) \binom{\lfloor \frac{a}{d} \rfloor - \lfloor \frac{a-1}{d} \rfloor}{\alpha}. \]
  Assume now that $A= \{a_1,a_2,\ldots,a_k\}$ and that the identities are true for $\{a_2,\ldots,a_k\}$. Without loss of generality
  we may assume that $a_1 < \sup A$. Then as to (a) we have
  \[
  \begin{split}
  f(\{a_1,\ldots,a_k\})
  &=
  f(\{a_2,\ldots,a_k\}) + \Phi(\{a_2,\ldots,a_k\},a_1) ) \\
  &=
  \sum_{d=1}^{\sup A}\mu(d) \left(2^{\sum_{i=2}^k
   (\lfloor \frac{a_i}{d}\rfloor - \lfloor \frac{a_i -1}{d}\rfloor)} -1\right)
   + \sum_{d|a_1} \mu(d)2^{\sum_{i=2}^k
   (\lfloor \frac{a_i}{d}\rfloor - \lfloor \frac{a_i -1}{d}\rfloor)} \\
  &=
  \sum_{d| a_1}\mu(d)
    \left(2^{\sum_{i=2}^k (\lfloor \frac{a_i}{d}\rfloor - \lfloor \frac{a_i -1}{d}\rfloor)} -1\right)
    +
    \sum_{\substack {d=1 \\ d \nmid a_1}}^{\sup A}\mu(d)
    \left(2^{\sum_{i=2}^k (\lfloor \frac{a_i}{d}\rfloor - \lfloor \frac{a_i -1}{d}\rfloor)} -1\right)
    + \\
   &  \quad\quad
   \sum_{d|a_1} \mu(d)2^{\sum_{i=2}^k (\lfloor \frac{a_i}{d}\rfloor - \lfloor \frac{a_i -1}{d}\rfloor)} \\
  &=
  2 \sum_{d|a_1} \mu(d)2^{\sum_{i=2}^k (\lfloor \frac{a_i}{d}\rfloor - \lfloor \frac{a_i -1}{d}\rfloor)}
    - \sum_{d|a_1}\mu(d) +
    \sum_{\substack {d=1 \\ d \nmid a_1}}^{\sup A}\mu(d)
    \left(2^{\sum_{i=2}^k (\lfloor \frac{a_i}{d}\rfloor - \lfloor \frac{a_i -1}{d}\rfloor)} -1\right) \\
  &=
  \sum_{d|a_1}\mu(d) \left(
    2^{\sum_{i=1}^k (\lfloor \frac{a_i}{d}\rfloor - \lfloor \frac{a_i -1}{d}\rfloor)} -1\right)
    +
    \sum_{\substack {d=1 \\ d \nmid a_1}}^{\sup A}\mu(d)
    \left(2^{\sum_{i=1}^k (\lfloor \frac{a_i}{d}\rfloor - \lfloor \frac{a_i -1}{d}\rfloor)} -1\right) \\
  &=
    \sum_{d=1}^{\sup A}\mu(d) \left(2^{\sum_{i=1}^k
   (\lfloor \frac{a_i}{d}\rfloor - \lfloor \frac{a_i -1}{d}\rfloor)} -1\right).
  \end{split}
  \]
  As to (b) we find
   \[
   \begin{split}
   f_{\alpha}(\{a_1,\ldots,a_k\})
  &=
  f_{\alpha}(\{a_2,\ldots,a_k\}) + \Phi_{\alpha -1}(\{a_2,\ldots,a_k\},a_1 ) \\
  &=
  \sum_{d=1}^{\sup A}\mu(d) \binom{\sum_{i=2}^k (\lfloor \frac{a_i}{d}\rfloor - \lfloor \frac{a_i -1}{d}\rfloor)}{\alpha}
   +
  \sum_{d|a_1}\mu(d) \binom{\sum_{i=2}^k (\lfloor \frac{a_i}{d}\rfloor - \lfloor \frac{a_i -1}{d}\rfloor)}{\alpha -1} \\
  &=
  \sum_{d|a_1}\mu(d)\left( \binom{\sum_{i=2}^k (\lfloor \frac{a_i}{d}\rfloor - \lfloor \frac{a_i -1}{d}\rfloor)}{\alpha}
  + \binom{\sum_{i=2}^k (\lfloor \frac{a_i}{d}\rfloor - \lfloor \frac{a_i -1}{d}\rfloor)}{\alpha -1} \right) + \\
  & \quad\quad
   \sum_{\substack{d=1\\ d\nmid a_1}}^{a_k}\mu(d) \binom{\sum_{i=1}^k (\lfloor \frac{a_i}{d}\rfloor - \lfloor \frac{a_i -1}{d}\rfloor)}{\alpha} \\
  &=
  \sum_{d|a_1}\mu(d) \binom{1+ \sum_{i=2}^k (\lfloor \frac{a_i}{d}\rfloor - \lfloor \frac{a_i -1}{d}\rfloor)}{\alpha}
  +
  \sum_{\substack{d=1\\ d\nmid a_1}}^{\sup A}\mu(d) \binom{\sum_{i=2}^k (\lfloor \frac{a_i}{d}\rfloor - \lfloor \frac{a_i -1}{d}\rfloor)}{\alpha} \\
  &=
  \sum_{d|a_1}\mu(d) \binom{\sum_{i=1}^k (\lfloor \frac{a_i}{d}\rfloor - \lfloor \frac{a_i -1}{d}\rfloor)}{\alpha}
  +
  \sum_{\substack{d=1\\ d\nmid a_1}}^{\sup A}\mu(d) \binom{\sum_{i=1}^k (\lfloor \frac{a_i}{d}\rfloor - \lfloor \frac{a_i -1}{d}\rfloor)}{\alpha} \\
  &=
  \sum_{d=1}^{\sup A}\mu(d) \binom{\sum_{i=1}^k (\lfloor \frac{a_i}{d}\rfloor - \lfloor \frac{a_i -1}{d}\rfloor)}{\alpha}.
  \end{split}
  \]
   This completes the proof.
 \end{proof}
 \noindent
 Alternatively, we have the following formulas for $f(A)$ and $f_{\alpha}(A)$.
 \begin{theorem} \label{main3}
 Let $A= \{a_1,a_2,\ldots,a_k\}$, let $\tau$ be a permutation of $\{1,2,\ldots,k\}$, and let
 $A_{\tau(j)} =\{a_{\tau(1)},a_{\tau(2)},\ldots, a_{\tau(j)} \}$ for $j=1,2,\ldots,k$. Then
 \[ (a)\quad
  f(A)= \sum_{j=1}^k \sum_{ d|a_{\tau(j)} } \mu(d) 2^{v(A_{\tau(j-1)},d)},
 \]
 \[ (b) \quad
 f_{\alpha}(A) = \sum_{j=1}^k \sum_{ d|a_{\tau(j)} } \mu(d) \binom{v(A_{\tau(j-1)},d)}{\alpha -1}.
 \]
 \end{theorem}
 \begin{proof}
 For simplicity we assume that $\tau$ is the identity permutation.
 As to part (a) we have
 \[
  \begin{split}
  f(\{a_1,\ldots,a_k\})
  &=
  f(\{a_1,\ldots,a_{k-1}\}) + \Phi(\{a_1,\ldots,a_{k-1}\},a_k)  \\
  &=
  f(\{a_1\}) + \Phi(\{a_1\},a_2) + \ldots + \Phi(\{a_1,\ldots,a_{k-1}\},a_k)  \\
  &=
  \sum_{d|a_1} \mu(d) + \sum_{d|a_2}\mu(d)
   2^{\left \lfloor \frac{a_1}{d} \rfloor - \lfloor \frac{a_1 -1}{d} \right\rfloor} + \ldots +
   \sum_{d|a_k} \mu(d) 2^{\sum_{i=1}^{k-1}
   (\lfloor \frac{a_i}{d}\rfloor - \lfloor \frac{a_i -1}{d}\rfloor)} \\
  &=
   \sum_{j=1}^k \sum_{d|a_j} \mu(d) 2^{\sum_{i=1}^{j-1}
   (\lfloor \frac{a_i}{d}\rfloor - \lfloor \frac{a_i -1}{d}\rfloor)},
  \end{split}
  \]
  where the third formula follows from Theorem \ref{main1}.
  Part (b) follows similarly. This completes the proof.
 \end{proof}
 \section{Combinatorial identities involving Mertens function}
\noindent
 We now give some identities which involves Mertens function.
 \begin{theorem} \label{Mertens-m,n}
 Let $m$ and $n$ be positive integers such that $1 < m \leq n$. Then
 \[  \sum_{d=1}^n \mu(d)
  (2^{\lfloor\frac{n}{d} \rfloor - \lfloor \frac{n-1}{d} \rfloor + \lfloor\frac{m}{d} \rfloor - \lfloor \frac{m-1}{d} \rfloor} )=
  \begin{cases}
 M(n)\ \text{if $(m,n) >1$,} \\
 1 + M(n)\ \text{if $(m,n)=1$.}
 \end{cases}
 \]
 \end{theorem}
 \begin{proof}
 Apply Theorem \ref{main2} (a) to the set $\{m,n \}$.
 \end{proof}
 \begin{theorem} \label{Mertens-l,m,n}
 Let $l$, $m$, and $n$ be positive integers such that $1 < l< m \leq n$. Then
 \[  \sum_{d=1}^n \mu(d)
  (2^{\lfloor\frac{n}{d} \rfloor - \lfloor \frac{n-1}{d} \rfloor + \lfloor\frac{m}{d} \rfloor - \lfloor \frac{m-1}{d} \rfloor +
   \lfloor\frac{l}{d} \rfloor - \lfloor \frac{l-1}{d} \rfloor}) = \]
 \[
  \begin{cases}
 4+ M(n)\ \text{if $(l,m)=(l,n)=(m,n) =1$,} \\
 3+ M(n)\ \text{if exactly two pairs from $\{l,m,n\}$ are co-prime,}\\
 2+ M(n)\ \text{if exactly one pair from $\{l,m,n\}$ is co-prime,}\\
 1 + M(n)\ \text{if no pair from $\{l,m,n\}$ is co-prime and $(l,m,n)=1$,} \\
 M(n)\ \text{Otherwise.}
 \end{cases}
 \]
 \end{theorem}
 \begin{proof}
 Apply Theorem \ref{main2} (a) to the set $\{l,m,n \}$.
 \end{proof}
 \begin{theorem} \label{Mertens-for-composite}
 We have
 \[
 M(\sup A) \leq \sum_{d=1}^{\sup A} \mu(d) 2^{v(A,d)},
 \]
 and equality occurs if and only if $\gcd(A) > 1$.
 \end{theorem}
 \begin{proof}
 The inequality follows
 Since by Theorem \ref{main2} (a) the identity
 \[
 \sum_{d=1}^{\sup A} \mu(d) 2^{v(b A,d)} - M(\sup A)
 \]
 counts the number of relatively prime subsets of $A$. Moreover, we clearly have
 that $\gcd (A) > 1$ if and only if $f(A) = 0$, which by Theorem \ref{main2} (a)
 means that
 \[
 \sum_{d=1}^{\sup A} \mu(d) 2^{v(b A,d)} = M(\sup A),
 \]
 as desired.
 \end{proof}
 \section{$\alpha$-relatively prime and $\alpha$-free relatively prime sets}
 \noindent
 Motivated by the work of Cameron and Erd\H{o}s in \cite{Cameron-Erdos} and the work of
 Calkin and Granville in \cite{Calkin-Granville} we introduce the following notions.
 \begin{definition}
 We say that the set $A$ is:
 \begin{itemize}
 \item \emph{$\alpha$-relatively prime} if every subset of $A$ of cardinality $\alpha$ is relatively prime,
 \item \emph{$\alpha$-relatively prime to $n$} if every subset of $A$ of cardinality $\alpha$ is relatively prime to $n$,
 \item \emph{$\alpha$-free relatively prime} if no subset of $A$ of cardinality $\alpha$ is relatively prime,
 \item \emph{$\alpha$-free relatively prime to $n$} if no subset of $A$ of cardinality $\alpha$ is relatively prime to $n$.
 \end{itemize}
 \end{definition}
\noindent
 We have the following algebraic characterizations.
 \begin{theorem} \label{alpha-relprime-characterization}
 (a) The set $A$ is $\alpha$-relatively prime if and only if $\binom{|A|}{\alpha} = f_{\alpha}(A)=0$. \\
 (b) The set $A$ is $\alpha$-relatively prime to $n$ if and only if $\binom{|A|}{\alpha} = \Phi_{\alpha}(A,n)$ . \\
 (c) The set $A$ is $\alpha$-free relatively prime if and only if $f_{\alpha}(A)=0$. \\
 (d) The set $A$ is $\alpha$-free relatively prime to $n$ if and only if $\Phi_{\alpha}(A,n) = 0$.
 \end{theorem}
 \begin{proof}
 Immediate from the definitions.
 \end{proof}
 \begin{corollary} \label{alpha-beta-relprime}
 Let $\beta$ be a positive integer such that $\alpha \leq \beta \leq |A|$. Then
 we have the following implications. \\
 \[
  \begin{split}
 \text{(a)  If\ } &
 \binom{|A|}{\alpha}=\sum_{d=1}^{\sup A}\mu(d)\binom{v(A,d)}{\alpha},\ \text{then\ }
  \binom{|A|}{\beta}=\sum_{d=1}^{\sup A}\mu(d)\binom{v(A,d)}{\beta}. \\
 \text{(b)  If\ } &
  \binom{|A|}{\alpha}=\sum_{d|n} \mu(d)\binom{v(A,d)}{\alpha},\ \text{then\ }
  \binom{|A|}{\beta}=\sum_{d|n}\mu(d)\binom{v(A,d)}{\beta}. \\
 \text{(c)  If\ } &
 \sum_{d=1}^{\sup A}\mu(d)\binom{v(A,d)}{\beta} =0,\ \text{then\ }
  \sum_{d=1}^{\sup A}\mu(d)\binom{v(A,d)}{\alpha} =0. \\
 \text{(c)  If\ } &
 \sum_{d|n}\mu(d)\binom{v(A,d)}{\beta} =0,\ \text{then\ }
  \sum_{d|n}\mu(d)\binom{v(A,d)}{\alpha} =0.
  \end{split}
  \]
 \end{corollary}
 \begin{proof}
 (a) Clearly if $A$ is $\alpha$-relatively prime and $\alpha \leq \beta \leq |A|$, then $A$ is $\beta$-relatively prime.
  Combining this fact with Theorem \ref{alpha-relprime-characterization} (a) and Theorem \ref{main2} (b) gives the desired implication. \\
  (b) Similarly, if $A$ is $\alpha$-relatively prime to $n$ and $\alpha \leq \beta \leq |A|$, then $A$ is $\beta$-relatively prime to $n$.
  Now combine this fact with Theorem \ref{alpha-relprime-characterization} (b) and Theorem \ref{main1} (b) to get the desired implication. \\
  (c) As to part (c), we use the fact that $A$ is $\alpha$-free relatively prime whenever $A$ is $\beta$-free relatively prime and
  $\alpha \leq \beta$. \\
  (d) As to part (d), we use the fact that $A$ is $\alpha$-free relatively prime to $n$ whenever $A$ is $\beta$-free relatively prime to $n$
  and $\alpha \leq \beta$. 
 \end{proof}
 Using similar ideas we have:
 \begin{corollary} \label{subset-relprime}
 Let $B$ be a nonempty subset of $A$ such that $\alpha \leq |B|$. Then we have the following implications. \\
 \[
 \begin{split}
 \text{(a)  If\ } &
 \binom{|A|}{\alpha}=\sum_{d=1}^{\sup A}\mu(d)\binom{v(A,d)}{\alpha},\ \text{then\ }
 \binom{|B|}{\alpha}=\sum_{d=1}^{\sup B}\mu(d)\binom{v(B,d)}{\alpha}. \\
 \text{(b)  If\ } &
  \binom{|A|}{\alpha}=\sum_{d|n} \mu(d)\binom{v(A,d)}{\alpha},\ \text{then\ }
  \binom{|B|}{\alpha}=\sum_{d|n}\mu(d)\binom{v(B,d)}{\alpha}. \\
 \text{(c) If\ } &
 \sum_{d=1}^{\sup A}\mu(d)\binom{v(A,d)}{\alpha} = 0,\ \text{then\ }
  \sum_{d=1}^{\sup B}\mu(d)\binom{v(B,d)}{\alpha} = 0. \\
 \text{(d) If\ } &
 \sum_{d|n} \mu(d)\binom{v(A,d)}{\alpha} = 0,\ \text{then\ }
  \sum_{d|n} \mu(d)\binom{v(B,d)}{\alpha} = 0.
 \end{split}
 \]
 \end{corollary}
 %
 Note that $A$ is $2$-relatively prime if and only if the elements of $A$ are pairwise co-prime.
 We have the following two characterizations for $2$-relatively prime sets.
 \begin{theorem} \label{two-relprime}
 Let $A=\{a_1,a_2,\ldots,a_k \}$, let $\tau$ be a permutation of $\{1,2,\ldots,k\}$, and let
 $A_{\tau(j)} =\{a_{\tau(1)},a_{\tau(2)},\ldots, a_{\tau(j)} \}$ for $j=1,2,\ldots,k$. Then the following are equivalent: \\
 (a)\quad The set $A$ is $2$-relatively prime. \\
 (b) \[ 2 |A| -1 = \left( 1 + 4 \sum_{d=1}^{\sup A} \mu(d) v(A,d) (v(A,d)-1) \right)^{1/2}. \]
 (c) \[ 2 |A| -1 = \left( 1 + 8 \sum_{j=1}^k \sum_{d|a_{\tau(j)}} \mu(d) v(A_{\tau(j-1)},d) \right)^{1/2}. \]
 \end{theorem}
 \begin{proof}
 As before we assume that $\tau$ is the identity permutation. By Theorem \ref{alpha-relprime-characterization}, Theorem \ref{main2}, and Theorem \ref{main3},
 the set $A$ is $2$-relatively prime means that
 \[ \binom{|A|}{2} = \sum_{d=1}^{\sup A}\mu(d)\binom{v(A,d)}{2} \]
 or equivalently
 \[ \binom{|A|}{2} = \sum_{j=1}^k \sum_{d|a_j} \mu(d)\binom{v(A_{j-1},d)}{1}. \]
 The former identity is equivalent to the quadratic equation
 \[ |A|^2-|A| - \sum_{d=1}^{\sup A} \mu(d) v(A,d)(v (A,d)-1) = 0 \]
 which means that
 \[ |A| = \frac{1}{2} \left(1+ \left(1 + 4 \sum_{d=1}^{\sup A} \mu(d) v(A,d) (v(A,d)-1) \right)^{1/2} \right), \]
 showing the equivalence of (a) and (b).
 The latter identity is equivalent to the quadratic equation
 \[ |A|^2-|A| - 2 \sum_{j=1}^k \sum_{d|a_j} \mu(d) v(A_{j-1},d) = 0 \]
 which is equivalent to
 \[ |A| = \frac{1}{2}\left(1 + \left(1+8 \sum_{j=1}^k \sum_{d|a_j} \mu(d) v(A_{j-1},d) \right)^{1/2} \right), \]
 showing the equivalence of (a) and (c). This completes the proof.
 \end{proof}
 \begin{corollary} \label{cor:pairwise-coprime}
 The number of nonempty subsets of $A$ whose elements are pairwise co-prime is
 \[ \# \left\{B\subseteq A:\ B\not=\emptyset\ \text{and\ }
   2|B| -1 = \left( 1 + 4 \sum_{d=1}^{\sup B} \mu(d) v(B,d) (v(B,d)-1) \right)^{1/2} \right\}.
 \]
 \end{corollary}
 Similar to Theorem \ref{two-relprime} we have:
 \begin{theorem} \label{two-free-relprime}
 Let $A=\{a_1,a_2,\ldots,a_k \}$, let $\tau$ be a permutation of $\{1,2,\ldots,k\}$, and let
 $A_{\tau(j)} =\{a_{\tau(1)},a_{\tau(2)},\ldots, a_{\tau(j)} \}$ for $j=1,2,\ldots,k$. Then the following are equivalent: \\
 (a)\quad The set $A$ is $2$-free relatively prime. \\
 (b)\quad \[ \sum_{d=1}^{\sup A} \mu(d) v(A,d) (v(A,d)-1) = 0.\]
 (c) \quad \[ \sum_{j=1}^k \sum_{d|a_{\tau(j)}} \mu(d) v(A_{\tau(j-1)},d) = 0. \]
 \end{theorem}
 \begin{corollary} \label{cor:pairwise-free-coprime}
 The number of nonempty subsets of $A$ which are co-prime free is
 \[ \# \{B\subseteq A:\ B\not=\emptyset\ \text{and\ }
   \sum_{d=1}^{\sup B} \mu(d) v(B,d) (v(B,d)-1)=0 \}.
 \]
 \end{corollary}
 %

\end{document}